\documentclass[english,twoside,11pt,fleqn]{article} 
\evensidemargin=\oddsidemargin
\evensidemargin=\oddsidemargin
\usepackage[dvips]{graphicx}
\usepackage[cp1251]{inputenc}           
\usepackage[english,russian]{babel}             
\usepackage{amsmath,amsthm,amssymb,newlfont}
\usepackage{cite}                       

\usepackage[hyper]{amsbib}              

\renewcommand{\Re}{\mathop{\mathrm{Re}}\nolimits}

\usepackage{array,hhline}
\usepackage{caption}
\newcounter{MyTable}

\DeclareCaptionFormat{MyFormat}{#1#2#3}
\DeclareCaptionLabelSeparator{MySeparator}{.\\}
\newenvironment*{MyTable}[1]              
{                                         
  \begin{table}[!h]
    \extrarowheight=2 pt
    \tabcolsep=1 mm
    \stepcounter{MyTable}
    \caption{#1}

    \begin{center}
}
{
    \end{center}
  \end{table}
}

\begin{document}
\title {\bf A non-homogeneous method of third order for additive stiff systems of ordinary differential equations
\footnote{The work was supported by Russian Fond of Fundamental Research (project № 05-01-00579-a).}}
\date{}
\author{ Evgeny Novikov, Anton Tuzov}
\maketitle

\renewcommand{\abstractname}{Abstract}
\renewcommand{\tablename}{Table}
\renewcommand{\refname}{References}

\begin{abstract}
  In this paper we construct a third order method for solving additively split autonomous stiff systems
  of ordinary differential equations.
  The constructed additive method is L-stable with respect to the implicit part and
  allows to use an arbitrary approximation of the Jacobian matrix.
In opposite to our previous paper~\cite{ASODE}, the fourth stage is explicit.
So, the constructed method also has a good stability properties because of $L$-stability of the intermediate numerical formulas
in the fourth stage, but has a lower computational costs per step.
  Automatic stepsize selection based on local error and stability control are performed.
  The estimations for error and stability control have been obtained without significant additional computational costs.
  Numerical experiments show reliability and efficiency of the implemented integration algorithm.
\end{abstract}

\section{Introduction}
  Spatial discretization of continuum mechanics problems in partial differential equations
  by finite difference or finite element methods results in the Cauchy problem for the system of ordinary differential equations
  with an additively split right hand side function of the form:
  \begin{equation*} 
    y'=\varphi(t,y)+g(t,y), \quad y(t_0)=y_0, \quad t_0\leq t\leq t_k,
  \end{equation*}
  where~$\varphi(t,y)$ is a non-symmetrical term obtained from discretization of the first-order differential operator,
  $g(t,y)$~is a symmetrical term obtained from discretization of the second-order differential operator,
    $t$~is a independent variable.
  It is assumed that in the problem the vector-function~$g$ is a stiff term and $\varphi$ is a non-stiff term.

  Explicit Runge-Kutta methods have a bounded stability region and are suitable for non-stiff and mildly stiff problems only.
  $L$-stable methods are usually used for solving stiff problems.
  In the case of large-scale problems overall computational costs of $L$-stable methods
  are almost completely dominated by evaluations and inversions of the Jacobian matrix of a right hand side vector function.
  Overall computational costs can be significantly reduced by re-using the same Jacobian matrix over several integration steps (freezing the Jacobian).

  Freezing the Jacobian in  iterative methods has effect on  convergence speed of an iterative process only and
  doesn't lead to loss of accuracy. So, this approach is extensively used for implementation of these  methods.
  For Rosenbrock type methods and their modifications~\cite{Hairer} an approximation of the Jacobian matrix
  can lead to decreasing a consistency order.

  The system $y'=f(t,y)$ can be written in the form  $y'=[f(t,y)-By ]+By$,
  where~$B$ is some  approximation of the Jacobian matrix. Assume that stiffness is fully concentrated in the term $g(t,y)=By$,
  then the expression $\varphi(t,y)=f(t,y)-By$ can be interpreted as the non-stiff term~\cite{Cooper_2, Novikov_add2}.
  If the Cauchy problem is considered in the form $y'=[f(t,y)-By ]+By$ under construction of additive methods,
  then an arbitrary approximation of the Jacobian matrix can be used without decreasing the order of these methods.
  Additive methods constructed in this way allow both analytical and numerical computations of the Jacobian matrix.
  Note that the approximation of the Jacobian by a diagonal matrix is suitable for some mildly stiff problems.

  In this paper we construct a six-stage third order additive method that allows to use different kinds of approximation of the Jacobian matrix.
In opposite to our previous paper~\cite{ASODE}, the fourth stage  is explicit.
The constructed method also has a good stability properties because of $L$-stability of the intermediate numerical formulas
(with respect to the implicit part) in the fourth stage, but has a lower computational costs per step.
  The estimation of the error has been obtained on the base of an embedded additive method without any additional computational costs.
  The estimation of the maximum absolute eigenvalue of the Jacobian matrix has been obtained by a power method
  using only two additional computations of ~$\varphi (y)$. Hence, additional computational costs  will be negligible, especially for large-scale problems.
  These estimations are used for error and stability control correspondingly.
  Numerical experiments are performed showing the reliability and efficiency of the constructed method.

 \section{A numerical scheme for autonomous problems}
  Consider the Cauchy problem for an autonomous system of ordinary differential equations
  \begin{equation} \label{add_a_sys}
    y'=\varphi(y)+g(y), \quad y(t_0)=y_0, \quad t_0 \leq t \leq t_k,
  \end{equation}
  where~$y,\varphi$ and~$g$~are $N$-dimensional smooth vector-functions, $t$~is an independent variable.
  In the following, we assume that $g$ is a stiff term and $\varphi$ is a non-stiff term.
  Consider a six-stage numerical scheme for solving~\eqref{add_a_sys}:
  \begin{align} \label{scheme}
    y_{n+1}&=y_n+\sum\limits_{i=1}^6 {p_i k_i}\, , \notag\\
       k_1 &= h\varphi(y_n),                        \notag\\
    D_nk_2 &= h [\varphi(y_n)+g(y_n) ],  \notag\\
    D_nk_3 & = k_2,                               \notag\\
       k_4 &= h\varphi (y_n+\sum\limits_{j=1}^3 \beta_{4j}k_j )+h g (y_n+\sum\limits_{j=1}^3 \alpha_{4j}k_j ), \\
    D_nk_5 &= k_4+\gamma k_3,                     \notag\\
       k_6 &= h\varphi (y_n+\sum\limits_{j=1}^5 \beta_{6j}k_j ), \notag
  \end{align}
  where~$D_n=E-ahg'_n$, $g'_n={\displaystyle \partial g(y_n)}/{\displaystyle \partial y}$~is the Jacobian matrix of the function $g(y)$,
  $E$~is the identity matrix, $k_i$, {$1\leq i\leq 6$}, are  stages,
  $a, p_i, \alpha_{4j}, \beta_{4j}, \beta_{6j},\gamma$~are coefficients that have effect on  accuracy and stability properties
  of the scheme~\eqref{scheme}.

\section{The third order conditions}
  The Taylor series expansion of the approximate solution up to terms in~$h^3$ has the form
  \begin{align*} 
    y_{n+1}&= y_n+ \bigl(p_1+p_2+p_3+p_4+(\gamma+1)p_5+p_6  \bigr)h\varphi+ \bigl(p_2+p_3+p_4+(\gamma+1)p_5 \bigr)hg+\\
           &+ \bigl((\beta_{41}+\beta_{42}+\beta_{43})(p_4+p_5)+(\beta_{61}+\beta_{62}+\beta_{63}+\beta_{64}+(\gamma+1)\beta_{65} )p_6 \bigr)h^2\varphi'\varphi+\\
           &+ \bigl((\beta_{42}+\beta_{43})(p_4+p_5)+(\beta_{62}+\beta_{63}+\beta_{64}+(\gamma+1)\beta_{65} )p_6 \bigr)h^2\varphi'g+\\
           &+ \bigl[a \bigl(p_2+2p_3+(3\gamma+1)p_5 \bigr)+(\alpha_{41}+\alpha_{42}+\alpha_{43})(p_4+p_5) \bigr]h^2 g'\varphi+\\
           &+ \bigl[a \bigl(p_2+2p_3+(3\gamma+1)p_5 \bigr)+(\alpha_{42}+\alpha_{43})(p_4+p_5) \bigr]h^2 g'g+\\
           &+0.5 \bigl[(\beta_{41}\!+\!\beta_{42}\!+\!\beta_{43})^2(p_4\!+\!p_5)+ \bigl(\beta_{61}\!+\!\beta_{62}\!+\!\beta_{63}\!+\!\beta_{64}\!+\!(\gamma\!+\!1)\beta_{65} \bigr)^2 p_6 \bigr]h^3\varphi''\varphi^2+\\
           &+0.5 \bigl[(\beta_{42}+\beta_{43})^2(p_4+p_5)+ \bigl(\beta _{62}+\beta_{63}+\beta_{64}+(\gamma+1)\beta_{65} \bigr)^2 p_6 \bigr]h^3\varphi''g^2+\\
           &+ \bigl[(\beta_{42}+\beta_{43})(\beta_{41}+\beta_{42}+\beta_{43})(p_4+p_5)+  \bigl(\beta_{61}+\beta_{62}+\beta_{63}+\beta_{64}+\\
           &+(\gamma+1)\beta_{65} \bigr) \bigl(\beta_{62}+\beta_{63}+\beta_{64}+(\gamma+1)\beta_{65} \bigr)p_6 \bigr]h^3 \varphi''\varphi g +\\
           &+(\beta_{41}+\beta_{42}+\beta_{43})(\beta_{64}+\beta_{65})p_6 h^3{\varphi'}^2\varphi+(\beta_{42}+\beta_{43})(\beta_{64}+\beta_{65})p_6 h^3{\varphi'}^2g+\\
           &+\Bigl[a \Bigl( \bigl(\beta_{42}+2\beta_{43} \bigr)\bigl(p_4+p_5\bigr)+\bigl(\beta_{62}+2\beta_{63}+(3\gamma+1)\beta_{65} \bigr)p_6 \Bigr)+\\
           &+\bigl(\alpha_{41}+\alpha_{42}+\alpha_{43}\bigr)\bigl(\beta_{64}+\beta_{65}\bigr)p_6 \Bigr]h^3 \varphi'g'\varphi+\Bigl[a \Bigl( \bigl(\beta_{42}+2\beta_{43} \bigr)\bigl(p_4+p_5\bigr)+\\
           &+\bigl(\beta_{62}+2\beta_{63}+(3\gamma+1)\beta_{65} \bigr)p_6 \Bigr)+\bigl(\alpha_{42}+\alpha_{43}\bigr)\bigl(\beta_{64}+\beta_{65}\bigr)p_6 \Bigr]h^3 \varphi'g'g+\\
           &+0.5(\alpha_{41}+\alpha_{42}+\alpha_{43})^2(p_4+p_5)h^3 g''\varphi^2+0.5(\alpha_{42}+\alpha_{43})^2(p_4+p_5)h^3 g''g^2+\\
           &+(\alpha_{41}+\alpha_{42}+\alpha_{43})(\alpha_{42}+\alpha_{43})(p_4+p_5)h^3 g''\varphi g+a(\beta_{41}+\beta_{42}+\beta_{43})p_5 h^3 g'\varphi'\varphi+\\
           &+a(\beta_{42}+\beta_{43})p_5 h^3 g'\varphi'g+a \bigl[a \bigl(p_2+3p_3+(6\gamma+1)p_5 \bigr)+(\alpha_{42}+2\alpha_{43})p_4+\\
           &+(\alpha_{41}+2\alpha_{42}+3\alpha_{43})p_5 \bigr]h^3 g'^2 \varphi+a \bigl[a \bigl(p_2+3p_3+(6\gamma+1)p_5 \bigr)+(\alpha_{42}+2\alpha_{43})p_4+\\
           &+(2\alpha_{42}+3\alpha_{43})p_5 \bigr]h^3 g'^2 g+O(h^4).\\
  \end{align*}
  where the corresponding elementary differentials are evaluated at~$y_n$.

\noindent
  The Taylor series expansion of the exact solution up to third order terms is
  \begin{align} \label{y(t)_acc3}
    y(t_{n+1})&=y(t_n)+h(\varphi+g)+\frac{h^2}{2}(\varphi'\varphi+\varphi'g+g'\varphi+g'g)+\frac{h^3}{6}(\varphi''\varphi^2+ \notag\\
              &+\varphi''g^2+2\varphi''\varphi g+{\varphi'}^2\varphi+{\varphi'}^2 g+\varphi'g'\varphi+\varphi'g'g+g''\varphi^2+g''g^2+                  \\
              &+2g''\varphi g+g'\varphi'\varphi+g'\varphi'g+{g'}^2\varphi+{g'}^2 g)+O(h^4 ),                                                                                                             \notag
  \end{align}
  where the corresponding elementary differentials are evaluated at~$y(t_n)$.

  Comparing the successive terms in the Taylor series expansion of the approximate and the exact solutions up to third order terms
  under the assumption $y_n=y(t_n)$ we have the system of nonlinear algebraic equations.
  Its solving results in the relation $\beta_{41}=\alpha_{41}=\beta_{61}=0$ and
  the third order conditions of the scheme~\eqref{scheme} take form:
  \begin{align} \label{sys2_acc3}
    &p_2+p_3+p_4+(\gamma+1)p_5=1, \notag\\
    &a(\beta_{42}+\beta_{43})p_5= 1/6, \notag\\
    & (\beta_{42}+\beta_{43})(\beta_{64}+\beta_{65})p_6=1/6, \notag\\
    & (\beta_{42}+\beta_{43})(p_4+p_5)+[\beta_{62}+\beta_{63}+\beta_{64}+(\gamma+1)\beta_{65}]p_6=0.5,\\
    & (\beta_{42}+\beta_{43})^2(p_4+p_5)+[\beta_{62}+\beta_{63}+\beta_{64}+(\gamma+1)\beta_{65}]^2 p_6=1/3,\notag\\
    &a(\beta_{42}+2\beta_{43})(p_4+p_5)\!+\!  \bigl[a \bigl(\beta_{62}\!+\!2\beta_{63}\!+\!(3\gamma\!+\!1)\beta_{65} \bigr)\!+\!(\alpha_{42}\!+\!\alpha_{43})(\beta_{64}\!+\!\beta_{65})  \bigr]p_6=1/6, \notag\\
    & (\alpha_{42}+\alpha_{43})^2(p_4+p_5)=1/3,\notag\\
    &a\bigl(p_2+2p_3+(3\gamma+1)p_5 \bigr) +(            \alpha_{42}+\alpha_{43})(p_4+p_5)= 0.5, \notag\\
    &a \bigl[a\bigl(p_2+3p_3+(6\gamma+1)p_5 \bigr)+(\alpha_{42}+2\alpha_{43})p_4+(2\alpha_{42}+3\alpha_{43})p_5 \bigr]=1/6, \notag\\
    &\alpha_{41}=\beta_{41}=\beta _{61}=0, \quad p_1=-p_6. \notag
  \end{align}

\section{Stability analysis}
  The linear stability analysis of the additive scheme~\eqref{scheme} is based on the scalar model equation
  \begin{equation} \label{eq_st}
    y'=\lambda_1 y + \lambda_2 y, \quad y(0)=y_0, \quad t \ge 0,\ \Re( \lambda_1) \leq 0,\ \Re( \lambda_2) \leq 0,\ |\Re( \lambda_1)| \ll |\Re( \lambda_2)|,
  \end{equation}
  where the free parameters $\lambda_1,\ \lambda_2$ can be interpreted as some eigenvalues of the Jacobian matrices of
  the functions $\varphi$ (the non-stiff term) and $g$ (the stiff term) correspondingly.

  Application of the scheme~\eqref{scheme} for numerical solving the equation~\eqref{eq_st} yields
  \begin{equation*}
    y_{n+1}=R(x,z)y_n,
  \end{equation*}
  where $x=\lambda_1 h, \ z=\lambda_2 h$ and $R(x,z)$ is a stability function (its analytical expression is omitted here for brevity).

  The necessary condition of $L$-stability of the additive scheme~\eqref{scheme} with respect to the stiff term has the form:
  \begin{equation*}
    \lim\limits_{z \to -\infty} R(x,z)=0.
  \end{equation*}
  It is satisfied if the following conditions hold:
  \begin{align} \label{L3}
    &\alpha_{42}=a,\quad \beta_{42}=0,\notag\\
    &a^2(p_1+p_6)-a(\beta_{62}+\beta_{64})p_6+\alpha_{43}\beta_{64}p_6=0,\\
    &a^2-a(p_2+p_4)+\alpha_{43}p_4=0.\notag
  \end{align}

  {\bf Solving the system (\ref{sys2_acc3}),~(\ref{L3}). }
  In the following, we assume that
  $ \sum_{j = 1}^3 \alpha _{4j} =1$ , $ \beta_{62} = a$.
  The first relation ensures that $g (y_n+\sum\limits_{j=1}^3 \alpha_{4j}k_j )$ approximate $g(y(t_{n+1}))$ in the fourth stage
   and the other one improve stability properties of the intermediate numerical formula.

  Let us denote
  \begin{align*}
    &\beta_1=\beta_{64}+\beta_{65}, & &\beta_2=\beta_{62}+\beta_{63}+\beta_{64}+(\gamma+1)\beta_{65},\\
    &\beta_3=a \bigl(\beta_{62}+2\beta_{63}+(3\gamma+1)\beta_{65} \bigr)+\beta_{64}+\beta_{65}, & &u=a(\gamma-1)+1.
  \end{align*}

  Then after obvious simplifications the system~\eqref{sys2_acc3},~\eqref{L3} takes the form
  \begin{align} \label{full_sys2}
    &p_2+p_3+p_4+(\gamma+1)p_5=1, \notag\\
    &a\beta_{43}p_5=1/6, \notag\\
    &\beta_{43}(\beta_{64}+\beta_{65})p_6=1/6, \notag\\
    &\beta_{43}(p_4+p_5)+ \bigl(\beta_{62}+\beta_{63}+\beta_{64}+(\gamma+1)\beta_{65} \bigr)p_6=0.5, \notag\\
    &\beta_{43}^2(p_4+p_5)+ \bigl(\beta_{62}+\beta_{63}+\beta_{64}+(\gamma+1)\beta_{65} \bigr)^2 p_6=1/3, \notag\\
    &2a\beta_{43}(p_4+p_5)+ \Bigl(a \bigl(\beta_{62}+2\beta_{63}+(3\gamma+1)\beta_{65} \bigr)+\beta_{64}+\beta _{65} \Bigr)p_6=1/6,\\
    &p_4+p_5=1/3, \notag\\
    &ap_2+2ap_3+p_4+ \bigl(a(3\gamma+1)+1 \bigr)p_5=0.5, \notag\\
    &a \bigl(ap_2+3ap_3+(2-a)p_4+3(2a\gamma+1)p_5 \bigr)=1/6, \notag\\
    &a^2-ap_2+(1-2a)p_4=0, \notag\\
    &\alpha_{41}=\beta_{41}=\beta_{42}=\beta_{61}=0,\; \alpha_{42}=\beta_{62}=a,\; \alpha_{43}=1-a,\;\beta_{64}=a^2/(1-2a),\; p_1 =-p_6. \notag
  \end{align}


  Multiplying the first equation of~\eqref{full_sys2}  by ~$2$ and subtracting the result from the eight one we obtain $-ap_2+(1-2a)p_4+up_5=0.5-2a$.
  It follows from here and the tenth equation of~\eqref{full_sys2} that
  \begin{equation} \label{eq1}
    p_5=0.5(2a^2-4a+1)/u.
  \end{equation}

  We shall try to obtain an equation for $a$.
  For this purpose we divide the ninth equation by $a$ and subtract the eighth one from the result. As the result we obtain:
  $ap_3=(a-1)p_4+(a-3a\gamma-2)p_5+(1-3a)/(6a)$. Substituting this relation to the eighth equation we have
  $ - ap_2  = (2a - 1)p_4  - 3up_5  + (2 - 9a)/(6a)$.
  It follows from here and the tenth equation of~\eqref{full_sys2} that $(6a^3-9a+2)/(6a)-3up_5=0$.
  Substituting this relation to~\eqref{eq1} we obtain the following equation for $a$:
  \begin{equation} \label{eq_on_a}
    6a^3-18a^2+9a-1=0.
  \end{equation}
  Then from the second equation of~\eqref{full_sys2} we have
  \begin{equation} \label{eq3}
    \beta_{43}=1/(6ap_5),
  \end{equation}
  It follows from~\eqref{eq3} and the third equation of~\eqref{full_sys2} that
  \begin{equation} \label{eq4}
    \beta_1=ap_5/p_6.
  \end{equation}
  From~\eqref{eq4} and the notation $\beta_1=\beta_{64}+\beta_{65}$ we obtain
  \begin{equation} \label{eq5}
    \beta_{65}=ap_5/p_6-\beta_{64},
  \end{equation}
  from~\eqref{eq3} and the fourth and the seventh equations of~\eqref{full_sys2}, we have
  \begin{equation} \label{eq6}
    \beta_2=1/(2p_6)-1/(18ap_5 p_6).
  \end{equation}
  It follows from sixth and seventh equations of~\eqref{full_sys2} and~\eqref{eq3} that $\beta_3=1/(6p_6 )-1/(9p_5 p_6 )$.
  From~\eqref{eq5},~\eqref{eq6} and the relations $\beta_2=\beta_{62}+\beta_{63}+\beta_{64}+(\gamma+1)\beta _{65}$, $\beta_{62}=a$,
  we obtain
  \begin{equation} \label{eq7}
    \beta_{63}=0.5 \bigl(1-2a(\gamma+1)p_5 \bigr)/p_6-1/(18ap_5 p_6 )-a+\gamma \beta _{64}.
  \end{equation}
  It follows from $\beta_3=a \bigl(\beta_{62}+2\beta_{63}+(3\gamma+1)\beta _{65} \bigr)+\beta_{64}+\beta_{65}$, \eqref{eq5},~\eqref{eq7},
  the eleventh equations of~\eqref{full_sys2} that:
  $\Bigl( \bigl(a^2(\gamma-1)+a \bigr)p_5+a-1/6 \Bigr)/p_6 = \bigl(a^3(\gamma-1)+a^2 \bigr)/(1 - 2a)$.
  Using the notations  introduced above the last equation can be written in the form: $\bigl(  a u p_5+a-1/6 \bigr)/p_6 = a^2 u/(1 - 2a)$.
  Substituting~\eqref{eq1} to this equation, we have $(6a^3-12a^2+9a-1)/(6p_6 ) =  a^2 u/(1 - 2a)$.
  Next, using~\eqref{eq_on_a} we obtain the following expression for $p_6 $:
  \begin{equation} \label{eq8}
    p_6=(1-2a)/u.
  \end{equation}

  It follows from the obtained results,  the seventh equation of~\eqref{full_sys2} and~\eqref{eq1} that
  $p_4=\bigl(-6a^2+2a(\gamma+5)-1 \bigr)/(6u)$.
  Substituting the last relation to the tenth equation of~\eqref{full_sys2} we obtain the expression for~$p_2: \;$
  Подставляя последнее равенство в десятое уравнение~\eqref{full_sys2}, выразим~$p_2: \;$
  $p_2= \bigl(6a^3(\gamma+1)-4a^2(\gamma+5)+2a(\gamma+6)-1 \bigr)/(6au)$.
  It follows from here and from~\eqref{eq_on_a} that
  \begin{equation} \label{eq9}
    p_2=\bigl(2a^2 (7\gamma-1)-a(7\gamma-3)+\gamma \bigr)/(6au).
  \end{equation}
  Substituting the seventh equation of~\eqref{full_sys2}, \eqref{eq1} and~\eqref{eq9} to the first one we express~$p_3: \;$
  $p_3= \bigl(-6a^3\gamma+2a^2(\gamma-1)+a(4\gamma+1)-\gamma \bigr)/(6au)$.
  Now, using~\eqref{eq_on_a}, we have
  $p_3=\bigl(-2a^2(8\gamma+1)+a(13\gamma+1)-2\gamma \bigr)/(6au)$.
  From~\eqref{eq1},~\eqref{eq4} and~\eqref{eq8}, we have
  $\beta_1=0.5a(2a^2-4a+1)/(1-2a)$.
  It follows from here and from~\eqref{eq_on_a} that
  $\beta_1 =(6a^2-6a+1)/(6-12a)$.
  From the last relation, the eleventh equation of~\eqref{full_sys2} $\beta_{64}=a^2/(1-2a)$ and $\beta_1=\beta_{64}+\beta_{65}$
  we obtain $\beta_{65}=-(6a-1)/(6-12a)$.
  From~\eqref{eq1},~\eqref{eq_on_a} and~\eqref{eq3} we have
  \begin{equation} \label{eq10}
    \beta_{43}=u/(6a^2-6a + 1).
  \end{equation}
  Substituting $\beta_{64}=a^2/(1-2a)$,~\eqref{eq1} and~\eqref{eq8} to~\eqref{eq7} and using~\eqref{eq_on_a}, we obtain
  $\beta_{63}= \bigl(-2a^2\gamma^2+2(53a^2-35a+4)\gamma+142a^2-86a+9 \bigr)/\bigl(6(-18a^2+10a-1) \bigr)$.

  We shall try to obtain an equation for $\gamma$.
  From the fourth and the fifth equations of~\eqref{full_sys2} we obtain
  $\beta_2=2(1-\beta_{43}^2)/(3-2\beta_{43})$.
  It follows from here and from the fifth equations of~\eqref{full_sys2} that
  $p_6=(3-2\beta_{43})^2/(12-12\beta_{43}^2)$.
  Substituting~\eqref{eq10} to the last relation and using \eqref{eq_on_a}, we have
  \begin{equation} \label{eq11}
    p_6=\frac{994a^2-72a^3\gamma+4a^2\gamma(\gamma+16)-4a(\gamma+143)+67}{12\ \bigl(102a^2-a^2(\gamma-1)^2-2a(\gamma+29)+6 \bigr)}.
  \end{equation}
  Comparing~\eqref{eq8} with~\eqref{eq11}, and using~\eqref{eq_on_a}, we obtain the following equation for $\gamma$:
  \begin{equation} \label{equation_on_gamma}
    4a^3\gamma^3-12a(15a^2-10a+1)\gamma^2+3a(374a^2-228a+33)\gamma+4(813a^2-486a+175/3)=0. \\
  \end{equation}
\noindent
  Now, the coefficients of the $L$~-stable third order scheme~\eqref{scheme} can be  computed by the following formulas:
 \begin{align} \label{koef}
    \alpha_{41}&=\beta_{41}=\beta_{42}=\beta_{61}=0,                       & \alpha_{42}&=\beta_{62}=a, \quad \alpha_{43}=1-a,\notag\\
            p_1&=-(1-2a)/u,                                                &         p_2&=\bigl(2a^2 (7\gamma\!-\!1)\!-\!a(7\gamma\!-\!3)\!+\!\gamma \bigr)/(6au),\notag\\
            p_3&=\bigl(-2a^2(8\gamma+1)+a(13\gamma+1)-2\gamma \bigr)/(6au),&         p_4&=\bigl(-6a^2+2a(\gamma+5)-1 \bigr)/(6u),\notag\\
            p_5&=0.5(2a^2-4a+1)/u,                                         &         p_6&=(1-2a)/u,\notag\\
     \beta_{43}&=1/(6ap_5),                                                &     \beta_1&=ap_5/p_6,\\
        \beta_2&=2(1-\beta_{43}^2)/(3-2\beta_{43}),                        &  \beta_{64}&=a^2/(1-2a),\notag\\
     \beta_{65}&=\beta_1-\beta_{64},                                       &  \beta_{63}&=\beta_2-a-\beta_{64}-(\gamma+1)\beta_{65}.\notag
  \end{align}
  where $u=a(\gamma-1)+1$, the coefficients~$a$ and $\gamma$ is determined from the equations~\eqref{eq_on_a} and  \eqref{equation_on_gamma} correspondingly.

  The equation~\eqref{eq_on_a} has the following tree real roots:
  \begin{equation*}
  a_1\!=\!0.15898389998867, \quad a_2\!=\!0.43586652150845, \quad a_3\!=\!2.40514957850286.
  \end{equation*}
  The numerical experiments show that the root~$a_2$ is the most suitable.
  The equation~\eqref{equation_on_gamma}, in turn, has the following tree real roots under the condition $a=a_2$:
  \begin{equation*}
    \gamma_{2,\;1}=-4.51745281449726, \quad \gamma_{2,\;2}=-2.49646456973997, \quad \gamma_{2,\;3}=-1.02332630944762\;.
  \end{equation*}
  The numerical experiments show that the root~$\gamma_{2,\;1}$. is the most suitable.
  Therefore computational results will be given for $a=a_2$ and $\gamma=\gamma_{2,\;1}$.

  The corresponding coefficients of the $L$~-stable third order scheme~\eqref{scheme} take the form
  \begin{align} \label{koef1}
              a&=+0.43586652150846, &      \gamma&=-4.51745281449727, \notag\\
            p_1&=+0.09130146290929, &         p_2&=+0.49588787677190, &         p_3&=+0.75521774748189, \notag\\
            p_4&=+0.20395977226114, &         p_5&=+0.12937356107220, &         p_6&=-0.09130146290929, \notag\\
    \alpha_{41}&=0,                 & \alpha_{42}&=+0.43586652150846, & \alpha_{43}&=+0.56413347849154, \notag\\
     \beta_{41}&=0,                 &  \beta_{42}&=0,                 &  \beta_{43}&=+2.95562753995095, \notag\\
     \beta_{61}&=0,                 &  \beta_{62}&=+0.43586652150846, &  \beta_{63}&=-3.98487214709651, \notag\\
     \beta_{64}&=+1.48112677684356, &  \beta_{65}&=-2.09874671679705\;.
  \end{align}

\section{Local error estimation}
  For the error estimation we construct the embedded method of second order of the form:
   \begin{align} \label{aux_scheme}
       y_{n+1,\;2}&=y_n+\sum\limits_{i=1}^3 {r_i k_i}+\sum\limits_{i=4}^6 {r_i \widetilde{k_i}}\;, \notag\\
               k_1&= h\varphi(y_n),                            \notag\\
            D_nk_2&= h \bigl(\varphi(y_n)+g(y_n) \bigr),      \notag\\
            D_nk_3& = k_2,                                         \\
   \widetilde{k_4}&= h\varphi (y_n+\sum\limits_{j=1}^3 \beta_{4j}k_j ), \notag\\
D_n\widetilde{k_5}&=\widetilde{k_4}+\gamma k_3, \notag\\
   \widetilde{k_6}&=h\varphi (y_n+\sum\limits_{j=1}^3 \beta_{6j}k_j +\beta_{64}\widetilde{k_4}+\beta_{65}\widetilde{k_5}). \notag
  \end{align}
  where the coefficients~$r_i$, {$1\leq i\leq 6$}, should be determined, and parameters $a,\beta_{4j}, \beta_{6j},\gamma$
  are given by~\eqref{koef} or~\eqref{koef1}.
  Note that there is not $h g (y_n+\sum_{j=1}^3 \alpha_{4j}k_j )$ in the fourth stage as opposed to~\eqref{scheme}.

  The Taylor series expansion of the approximate solution computed by the scheme~\eqref{aux_scheme} up to terms in~$h^2$ has the form
  \begin{align*} 
    y_{n+1,\;2}&=y_n+ \bigl(r_1+r_2+r_3+r_4+(\gamma+1)r_5+r_6 \bigr)h\varphi + (r_2+r_3+\gamma r_5)hg + \\
               & + a \bigl(r_2+2r_3+(3\gamma+1)r_5 \bigr)h^2 g'\varphi + a(r_2+2r_3+3\gamma r_5)h^2 g'g + \\
               & + \Bigl((\beta_{41}+\beta_{42}+\beta_{43})(r_4+r_5)+ \bigl(\beta_{61}+\beta_{62}+\beta_{63}+\beta_{64}+(\gamma+1)\beta_{65} \bigr)r_6  \Bigr)h^2 \varphi'\varphi + \\
               & + \bigl((\beta_{42}+\beta_{43})(r_4+r_5)+(\beta_{62}+\beta_{63}+\gamma\beta_{65})r_6 \bigr)h^2 \varphi'g + O(h^3 ).\\
  \end{align*}
  where the elementary differentials are evaluated at~$y_n$.
  Comparing successive terms in the Taylor series expansion of the approximate and the exact solutions up to second order terms
  under the assumption $y_n=y(t_n)$ we obtain the second order conditions of the scheme~\eqref{aux_scheme}:
  \begin{align} \label{sys_acc2}
    1.\ \ &r_1+r_2+r_3+r_4+(\gamma+1)r_5+r_6=1, \notag\\
    2.\ \ &r_2+r_3+\gamma r_5=1, \notag\\
    3.\ \ &a \bigl(r_2+2r_3+(3\gamma+1)r_5 \bigr)=0.5, \notag\\
    4.\ \ &a(r_2+2r_3+3\gamma r_5)=0.5, \\
    5.\ \ &(\beta_{41}+\beta_{42}+\beta_{43})(r_4+r_5)+ \bigl(\beta_{61}+\beta_{62}+\beta_{63}+\beta_{64}+(\gamma+1)\beta_{65} \bigr)r_6=0.5, \notag\\
    6.\ \ &(\beta_{42}+\beta_{43})(r_4+r_5)+(\beta_{62}+\beta_{63}+\gamma \beta_{65})r_6=0.5. \notag
  \end{align}
  Subtracting the third equation of~\eqref{sys_acc2} from the fourth one and using $a \ne 0$ we obtain $r_5=0$.
  Subtracting the sixth equation from the fifth one and using the eleventh one of ~\eqref{full_sys2} and  $\beta_{64}+\beta_{65} \ne 0$
  we have~$r_6=0$.
  It follows from the second and the fourth equations of~\eqref{sys_acc2} that
  $r_2=0.5(4a-1)/a$, $r_3=0.5(1-2a)/a$.
  Now, from the sixth and the first equation of~\eqref{sys_acc2} we have
  $r_4=0.5/ \beta_{43}$, $r_1=- 0.5/ \beta_{43}$.

  As the result we have all the coefficients of the $L$~stable embedded method~\eqref{aux_scheme} of second order.
  For the coefficients~\eqref{koef1} we obtain
  \begin{align*}
    r_1&=-0.16916881211910, & r_2&=+0.85285981986048   & r_3&=+0.14714018013952,\\
    r_4&=+0.16916881211910, & r_5&=0.                  & r_6&=0.
  \end{align*}

  The embedded method~\eqref{aux_scheme} doesn't require any additional computations of right hand side,
  evaluations and inversions of the Jacobian matrix,  because of $r_6=0$ and ${r_5}=0$.

  Let us denote the error estimation by
  \begin{equation*}
    err_n=\max\limits_{1 \leq i \leq N}\dfrac{ |y_n^{i}-y_{n,\;2}^{i}|}{Atol_i+Rtol_i|y_n^{i}|},
  \end{equation*}
  where $Atol_i$ and $Rtol_i$ are the desired tolerances prescribed by the user.
  If $ err_n  \leq 1,$ then the computed step is accepted, else the step is rejected and computations are repeated.
  When $Rtol_i=0$, the absolute error is controlled on the $i$-th component of the solution with the desired tolerance $Atol_i$.
  If $Atol_i=0$ then the relative error is controlled on the $i$-th component with the tolerance $Rtol_i$.

\newpage
\section{Stability control and stepsize selection}
  In the additive method~\eqref{scheme} for solving~\eqref{add_a_sys}
  the non-stiff term~$\varphi$ is treated by the tree-stage explicit Runge Kutta method (the explicit part),
  and the stiff term~$g$ is treated by the $L$-stable $(4,2)$-method~\cite{Novikov_freezing, Nov_1, Novikov_m_kk1} (the implicit part).
  In the general case there is no guarantee that the function $\varphi (y)=f(y)-By$ is the non-stiff term
  in reducing $y'=f(y)$ to $y'=[f(y)-By]+By$.
  If some stiffness is in $\varphi (y)=f(y)-By$ (i.e. stiffness leakage phenomenon occurs)
  then the additional stability control of the explicit part of the scheme~\eqref{scheme}
  can increase efficiency of computations for many problems.
  In some cases it has no a significant effect on the efficiency of the integration algorithm because of
  the good stability properties of the scheme~\eqref{scheme}.
  Therefore the choice of using or not using  the additional stability control of the explicit part is given to the end-user.

  We perform the stability control of the explicit part of the scheme~\eqref{scheme} by analogy with~\cite{Novikov_monograph}.
  For additive methods in opposite to explicit Runge Kutta methods it isn't possible to use previously computed stages
  because of  peculiarity of the problem~\eqref{add_a_sys}.
  Therefore instead of using the stages $k_i, \ 1 \leq i \leq 6,$ of~\eqref{scheme}
  we consider the additional stages~$d_1$, $d_2$ of the form:
  \begin{equation*}
    d_1=h\varphi (y_n+\alpha_{21}k_1), \quad d_2=h\varphi (y_n+\alpha_{31}k_1+\alpha_{32}d_1).
  \end{equation*}
  Denote $\varphi (y)=Ay+b$, where~$A$ and~$b$ are matrix and vector with constant coefficients correspondingly, then we have
  \begin{equation*}
    k_1=h(Ay_n+b), \quad d_1=k_1+\alpha_{21}hAk_1, \quad d_2=k_1+(\alpha_{31}+\alpha_{32})hAk_1+\alpha_{21}\alpha_{32}h^2A^2k_1.
  \end{equation*}
  Assuming $\alpha_{21}=\alpha_{31}+\alpha_{32}$ we obtain
  \begin{equation*}
    d_2-d_1=\alpha_{21}\alpha_{32}h^2A^2k_1, \quad d_1-k_1=\alpha_{21}hAk_1.
  \end{equation*}
  The maximum absolute eigenvalue~$v_n=h|\lambda_{n\ max}|$ of the matrix~$hA$ can be approximated using the power method by the following formula:
  \begin{equation*}
    v_n=|\alpha_{32}^{-1}|\max\limits_{1 \leq i \leq N}\dfrac{|d_2^i-d_1^i|}{|d_1^i-k_1^i|},
  \end{equation*}
  then the stability control can be made by  $v_n  \leq 2,$ where number 2 is an approximate length of the stability interval of
  the tree-stage explicit Runge Kutta method.

  In the general case this estimation is quite crude because of small number of iterations of the power method
  and the nonlinearity of the function~$\varphi (y)$. Therefore the stability control is used for limiting the stepsize growing only.

  Let the approximate solution ~$y_n$ is computed with the stepsize~$h_n$. For the stepsize selection we use
  $err_n = O(h_n^3)$. The stepsize~$h_{acc}$ predicted by accuracy we compute by the formula: $h_{acc}=q_1h_n$,
  where~$q_1$ is a root of the equation $q_1^3 err_n=1.$
  In view of $v_n=O(h_n)$, the stepsize~$h_{st}$ predicted by stability is computed by  $h_{st}=q_2h_n$,
  where~$q_2$ is a root of the equation $q_2v_n=2$. Then the stepsize~$h_{n+1}$ predicted by accuracy and stability is selected by the formula:
  \begin{equation*}
    h_{n+1}=\max[h_n, \min(h_{acc},h_{st})].
  \end{equation*}

  The stability control of the explicit part of the scheme~\eqref{scheme} requires, at each integration step,
  two additional computations of~$\varphi (y)$.  These computational costs are negligible for large-scale problems,
  but if you are sure that all stiffness is in~$g(y)$ then you can take off stability control to save computational costs.

\newpage
\section{Numerical experiments}
  Further, the numerical code based on the additive method~\eqref{scheme} (with error and stability control as well as with diagonal Jacobian approximation)
  is called ASODE3 (the Additive Solver of Ordinary Differential Equations).

  The test problems given below have been reduced to the form $y'=(f(y)-By)+By$.
  All numerical computations have been performed in double precision arithmetic on IBM~PC~Athlon(tm) XP~2000+
  with the desired tolerances of the error $Atol=Rtol=Tol=10^{-m},\ m=2,4$.
  The scheme~\eqref{scheme} is of third order, therefore it is unreasonable to do numerical computations  with higher tolerance.

  The following four test examples are considered:

  Example~1~\cite{Enright}.
  \begin{align} \label{f28}
    y'_1&=-0.013y_1-1000y_1y_3,             \notag\\
    y'_2&=-2500y_2y_3,                      \\
    y'_3&=-0.013y_1-1000y_1y_3-2500y_2y_3,  \notag\\
        &\ t \in [0,50], \quad y_1(0)=1, \quad y_2(0)=1, \quad y_3(0)=0, \quad h_0=2.9\cdot 10^{-4}. \notag
    \end{align}

  Example~2~\cite{Gear}.
  \begin{align} \label{f29}
    y'_1&=77.27(y_2-y_1y_2+y_1-8.375 \cdot 10^{-6}y_1^2), \notag\\
    y'_2&=(-y_2-y_1y_2+y_3)/77.27,\\
    y'_3&=0.161(y_1-y_3),\notag\\
        &\ t \in  [0,300], \quad y_1(0)=4, \quad y_2(0)=1.1, \quad y_3(0)=4, \quad h_0=2\cdot 10^{-3}. \notag
  \end{align}

  Example~3.
  \begin{eqnarray*}
    y'_1&=&-0.04y_1+0.01y_2y_3,\nonumber\\
    y'_2&=&400y_1-100y_2y_3-3000y_2^2, \nonumber\\
    y'_3&=&30y_2^2, \nonumber\\
    t &\in & [0,40], \quad y_1(0)=1, \quad y_2(0)=y_3(0)=0, \quad h_0=10^{-5}. \nonumber
  \end{eqnarray*}

  Example~4.
  \begin{eqnarray*}
    y'_1&=&y_3-100y_1y_2,\nonumber\\
    y'_2&=&y_3+2y_4-100y_1y_2-2\cdot 10^4y_2^2,\nonumber\\
    y'_3&=&-y_3+100y_1y_2, \nonumber\\
    y'_4&=&-y_4+10^4y_2^2,\nonumber\\
    t &\in & [0,20], \quad y_1(0)=y_2(0)=1, \quad y_3(0)=y_4(0)=0, \quad h_0=2.5\cdot 10^{-5}. \nonumber
  \end{eqnarray*}

  The approximation of the Jacobian by a diagonal matrix is used when solving the test problems by ASODE3.
  For the first test problem the diagonal matrix $B$ with elements
  $b_{11}=-0.013-1000y_3$, $b_{22}=-2500y_3$, $b_{33}=-1000y_1-2500y_2$ are used.
  In the case of diagonal Jacobian approximation computational costs of additive methods
  are dominated by the number of right hand side function evaluations.
  So, computational costs of~\eqref{scheme} per integration step are comparable to ones of explicit methods.
  Hence, ASODE3 is compared with the following numerical codes based on well-known explicit Runge Kutta methods:\\
  \vspace{-2\baselineskip}
  \begin{tabbing}
    RKM4 \= -- 5-stage \; \= Merson  method of order 4~\cite{Merson},\\
    RKF5 \> -- 6-stage    \> Felberg method of order 5~\cite{Fehlberg},\\
    RKF7 \> -- 13-stage   \> Felberg method of order 7~\cite{Fehlberg},\\
    DP8  \> -- 13-stage   \> Dormand and Prince method method of order 8~\cite{Dormand},\\
    and less well-known Runge-Kutta type method: \\
    RKN2 \> -- 2-stage    \> method of order 2~\cite{Novikov_DSc}.\\
  \end{tabbing}
  \vspace{-2\baselineskip}

  The overall computational costs  (measured by the number of right hand side function evaluations over the integration interval)
  are given in the table

\begin{MyTable}{Computational costs of RKM4, RKF5, RKF7, DP8, RKN2, ASODE3 with stability control.}
\label{tab_numerical_experiments}
\begin{tabular}{|c|c|r|r|r|r|r|r|r|r|}
\hline
\multicolumn{1}{|c|}{\;\; \No \;\;}  &  \multicolumn{1}{c|}{\; Tol \;}  & \multicolumn{1}{c|}{RKM4}  & \multicolumn{1}{c|}{RKF5} & \multicolumn{1}{c|}{RKF7}  & \multicolumn{1}{c| }{DP8} & \multicolumn{1}{c|}{RKN2}  & \multicolumn{1}{c|} {ASODE3}\\
\hline
 1   & $10^{-2}$ & 401 716    & 401 005    & 982 536    & 717 526    & 222 441      & 9 351       \\
      \cline{2-8}
      & $10^{-4}$ & 400 627    & 400 656    & 982 150    & 717 287    & 222 481     & 37 338    \\
\hline
  2   & $10^{-2}$ & 13 391 594 & 15 694 434 & 38 429 196 & 27 998 053 & 8 682 849   & 1 589   \\
      \cline{2-8}
      & $10^{-4}$ & 13 384 132 & 15 691 105 & 38 429 976 & 27 993 793 & 8 689 861   & 7 711 \\
\hline
  3   & $10^{-2}$ & 204 889    & 237 942    & 587 509    & 431 591    & 133 022     & 3 129     \\
      \cline{2-8}
      & $10^{-4}$ & 206 647    & 240 676    & 565 396    & 430 823    & 132 987     & 16 361    \\
\hline
  4   & $10^{-2}$ & 10 832     & 11 874     & 29 991     & 23 052     & 6 585       & 63 430     \\
      \cline{2-8}
      & $10^{-4}$ & 10 236     & 11 366     & 28 819     & 23 354     & 7 627       & 367 411   \\
\hline
\end{tabular}

\end{MyTable}

\section{Conclusions}

  In addition to continuum mechanics problems, the constructed additive method can be used for solving locally unstable problems.
  In this case~$\varphi (y)$ corresponds to eigenvalues of the Jacobian matrix with positive real parts.
  In opposite to $A$-stable methods, explicit Runge Kutta methods are unstable in almost the entire right half plane
  and therefore are more suitable for detecting the local unstable solutions.
  For many locally unstable problems it is also easy to split the right hand side into stiff and non-stiff terms from physical considerations.

  So, in this paper, we constructed the third order additive method that is $L$-stable with respect to the implicit part
  and allows to use an arbitrary approximation of the Jacobian matrix without loss of accuracy.
  Automatic stepsize selection based on local error and stability control are performed and
  the auxiliary formulas for doing this were obtained without significant additional computational costs.

  The aim of numerical computations was to test the reliability and efficiency of the implemented integration algorithm
  with error and stability control as well as with diagonal Jacobian approximation.
  They didn't aim at solving practical problems of continuum mechanics and locally unstable problems.
\noindent
  Numerical experiments show reliability and efficiency of the presented method.
  It follows from them that the method has good stability properties for solving mildly stiff problems
  and that the test problems turned out to be rather stiff for the explicit Runge-Kutta methods considered above.
  It is worth remarking that computational costs per step are comparable for both
  the additive method (with diagonal Jacobian approximation) and explicit ones.
  So, the implemented integration algorithm makes it possible  to expend the range of applicability
  of explicit Runge-Kutta methods towards more stiff problems.

\end{document}